\documentclass[12pt,a4paper]{article}

\usepackage[margin=2cm]{geometry}
\usepackage{amsmath,amssymb,latexsym,amsthm}

\theoremstyle{plain}
\newtheorem{proposition}{Proposition}[section]
\newtheorem{theorem}[proposition]{Theorem}

\newtheorem{lemma}[proposition]{Lemma}
\theoremstyle{definition}
\newtheorem{definition}[proposition]{Definition}
\newtheorem{remark}[proposition]{Remark}

\newcommand{\ip}[2]{\langle #1,#2 \rangle}
\newcommand{\mc}{\mathcal}

\newcommand{\G}{{\mathbb G}}
\newcommand{\vnten}{\overline\otimes}
\newcommand{\proten}{\widehat\otimes}

\begin{document}

\title{Completely positive multipliers of quantum groups}
\author{Matthew Daws}
\maketitle

\begin{abstract}
We show that any completely positive multiplier of the convolution algebra
of the dual of an operator algebraic quantum group $\G$ (either a locally
compact quantum group, or a quantum group coming from a modular or manageable
multiplicative unitary) is induced in a canonical fashion by a unitary
corepresentation of $\G$.  It follows that there is an order bijection between
the completely positive multipliers of $L^1(\G)$ and the positive functionals on
the universal quantum group $C_0^u(\G)$.  We provide a direct link between the
Junge, Neufang, Ruan representation result and the representing element of a
multiplier, and use this to show that their representation map is always
weak$^*$-weak$^*$-continuous.


\noindent\emph{MSC classification (2010):} 20G42, 22D10, 22D15, 43A22,
   46L07, 46L89, 81R50.
\end{abstract}

\section{Introduction}

Multipliers arise throughout the study of algebras in analysis, as a useful
tool for embedding a non-unital algebra into its ``largest'' unitisation.
In abstract harmonic analysis, the convolution algebra $L^1(G)$ of a locally
compact group $G$ is unital only when $G$ is discrete; Wendel's Theorem
tells us that the multiplier algebra of $L^1(G)$ is the measure algebra $M(G)$.
On the ``Fourier transform'' side, the Fourier algebra $A(G)$ is unital only
for compact $G$.  Here it seems most profitable to study the completely
bounded multipliers, as these are more tractable, have better functorial
properties, and can capture interesting geometric aspects of the group $G$,
see \cite{ch, dch, spronk} for example.

The convolution algebras of locally compact quantum groups $\G$ form a
generalisation of both $L^1(G)$ and $A(G)$, together with genuinely
``non-commutative'' examples, such as compact quantum groups.
Again, completely bounded multipliers are profitable to study,
and in particular, it was shown by Junge, Neufang and Ruan in \cite{jnr} that
there is a bijection between ``abstract'' completely bounded (left) multipliers
of the dual $L^1(\hat\G)$ (that is, right module maps on $L^1(\hat\G)$) and
``concrete'' multipliers-- elements of $L^\infty(\G)$ which ``multiply''
the image of $L^1(\hat\G)$ into itself.  This mirrors the Fourier algebra picture,
where there is a bijection between completely bounded multipliers of $A(G)$,
thought of as module maps on $A(G)$, and continuous functions $f$ on $G$ which
multiply the function algebra $A(G)$ into itself.

It was shown by De Canniere and Haagerup, \cite[Proposition~4.2]{dch}, that
completely positive multipliers of $A(G)$ biject with positive definite
functions on $G$.  Positive definite functions arise from unitary group
representations, equivalently, $*$-representations of the full group C$^*$-algebra.
In this light, we can also observe that positive multipliers of $L^1(G)$ arise
from positive measures in $M(G)$, or equivalently, $*$-representations of $C_0(G)$.
Similarly, the use of completely positive multipliers of quantum groups has
arisen in the study of various approximation properties associated to quantum
group algebras; see for example \cite{fre}, and in relation to the Haagerup
approximation property, \cite{b1,b2}.

The main objective of this paper is to establish the same result for quantum groups:
any completely positive multiplier of $L^1(\hat\G)$ comes from a unitary
corepresentation of $\G$, or equivalently, from a $*$-representation of the
universal quantum group $C_0^u(\hat\G)$ (to be thought of as generalising
the role of $C^*(G)$.)  We remark that such a result cannot in general be extended
to completely bounded multipliers, as unpublished work of Losert, extending
his paper \cite{los}, shows that $G$ is amenable if and only if the
Fourier-Stieltjes algebra $B(G)=C^*(G)^*$ agrees with the algebra of completely
bounded multipliers of $A(G)$.  (This was shown for various explicit examples
in \cite{ch, dch} and for all discrete groups in \cite{boz}.)

A key tool for us is the representation result of \cite{jnr},
where it is shown that a completely bounded multiplier of $L^1(\hat\G)$
extends to a normal, completely bounded map on $\mc B(L^2(\G))$ with certain
properties.  In the completely positive case, applying the Stinespring
construction to this map yields a Hilbert space, and it is this space which
will be the carrier for our unitary corepresentation.

Thus a completely bounded (left) multiplier gives rise to both a normal
map $\Phi$ on $\mc B(L^2(\G))$ and a ``representing element'' $a_0\in
L^\infty(\G)$.  We show that these are linked by the relation
\[ \big( \Phi(\theta_{\xi,\eta})\alpha \big| \beta \big)
= \ip{\Delta(a_0)}{\omega_{\alpha,\eta} \otimes
\omega_{\xi,\beta}^{\sharp*}}. \]
Here $\theta_{\xi,\eta}$ is the compact operator $\alpha\mapsto
(\alpha|\eta)\xi$, and $\sharp$ denotes the (in general unbounded)
$*$-operation on $L^1(\G)$.  See Proposition~\ref{prop:rep_gives_mult} for
a precise statement.  We use this result to show that the representation
map of \cite{jnr} is always weak$^*$-weak$^*$ continuous, which answers
in the affirmative the question asked before \cite[Theorem~4.7]{hnr2}.
Furthermore, this suggests a possible definition for
a ``positive definite function'' on a quantum group: namely, those $a_0$ such
that $\Phi$ can be chosen to be completely positive.  In joint work with
Salmi, \cite{dawssalmi}, we verify this conjecture-- namely, it is indeed the
case that if $\Phi$ can be chosen to be completely positive, then automatically
$\Phi$ must come from a completely positive \emph{multiplier}, and so $a_0$
comes from a corepresentation of $\G$.

The organisation of the paper is as follows.  In Section~\ref{sec:qg} we
quickly survey the Operator Algebraic quantum groups which we shall study,
and we prove some useful auxiliary results.  In this paper, all
our arguments work in the slightly more general setting of quantum groups coming
from manageable (or modular) multiplicative unitaries (that is, without assuming
invariant weights).  In Section~\ref{sec:mult} we
recall the construction of \cite{jnr}.  As the actual constructions involved,
and not just the results of that paper, are indispensable to us, we sketch
some of the proofs, and take the opportunity to show how the arguments can be
made to work for quantum groups coming from manageable multiplicative unitaries.
In Section~\ref{sec:corep-mult} we show how certain corepresentations of
$\G$ canonically give rise to completely bounded multipliers of $L^1(\hat\G)$.
We also show how unitary corepresentations, equivalently, $*$-representations
of the universal dual algebra $C_0^u(\hat\G)$, give rise to completely
positive multipliers of $L^1(\hat\G)$.  The following section is the main
part of the paper, and proves the converse-- any completely positive multiplier
of $L^1(\hat\G)$ comes from a unitary corepresentation of $\G$ in the manner
discussed in Section~\ref{sec:corep-mult}.  As mentioned above, we construct
the Hilbert space for our corepresentation from the Stinespring dilation
construction.  We can view this as a completion of $L^1(\G)$ with respect to
a certain inner-product; then the corepresentation $U$ is built by
actually forming $U^*$, which is linked to the anti-representation of $L^1(\G)$
given by right multiplication of $L^1(\G)$ on itself.  We then proceed to show
that $U$ is in fact unitary, and then that we can recover the original
multiplier from $U$.  In the final section, we explore,
as discussed above, the link between multipliers as maps on $L^1(\G)$
(or, via \cite{jnr}, as maps on $\mc B(L^2(\G))$) and their ``representing
elements'' in $L^\infty(\G)$.

A final word on notation.  Our Hilbert space inner products shall be linear
in the first variable, and we write $(\cdot|\cdot)$ for an inner product
(or more generally, a sesquilinear form).  We write $\ip{\cdot}{\cdot}$
for the bilinear pairing between a Banach space and its dual.
For a Hilbert space $H$, we write $\mc B(H)$ for the
algebra of all bounded operators on $H$, write $\mc B(H)_*$ for its predual
(the trace class operators) and write $\mc B_0(H)$ for the ideal of compact
operators.  For $\xi,\eta\in H$, we denote by $\omega_{\xi,\eta}$ the
normal functional in $\mc B(H)_*$, and by $\theta_{\xi,\eta}$ the rank one
operator in $\mc B_0(H)$, which are defined by
\[ \ip{T}{\omega_{\xi,\eta}} = (T\xi|\eta), \qquad
\theta_{\xi,\eta}(\gamma) = (\gamma|\eta) \xi
\qquad (T\in\mc B(H), \gamma\in H). \]

Given a normal map $T$ on a von Neumann algebra $M$, we write
$T_*$ for the pre-adjoint of $T$ acting on the predual $M_*$.  We write
$\otimes$ to mean a completed tensor product, either of Hilbert spaces, or
the minimal C$^*$-algebraic tensor product.  We write $\vnten$ for the
von Neumann algebraic tensor product, and $\odot$ for the purely algebraic
tensor product.  We write $\Sigma$ for the tensor swap map of Hilbert spaces,
say $\Sigma:H\otimes H\rightarrow H\otimes H; \xi\otimes\eta \mapsto
\eta\otimes\xi$.  

We use the basic theory of Operator Spaces without comment; see, for example,
\cite{er} for further details.

\medskip

\noindent\textbf{Acknowledgements:}
The author wishes to thank Michael Brannan for asking the initial question
which lead to this paper, and to thank Zhong-Jin Ruan for pointing out that
the argument at the end of Proposition~\ref{prop:mults_rep}
needed more justification.  The anonymous referee of an earlier version of
this paper made many helpful comments; in particular suggesting the argument of
Theorem~\ref{thm:gen_all_bh}.

\section{Operator algebraic quantum groups}\label{sec:qg}

In this paper, we shall be concerned with quantum groups in the operator
algebraic setting-- to be precise, either locally compact quantum groups,
in the Kustermans, Vaes sense \cite{kus1, kv, kvvn, vaesphd}, or
C$^*$-algebraic quantum groups built from manageable or modular
multiplicative unitaries, in the So{\l}tan, Woronowicz sense
\cite{sw2, sw, woro} (the latter generalising the former).  In fact, for
many of our results, we shall need remarkably little-- our main tool being
that ``invariants are constant'' (see below-- our inspiration here is
\cite[Section~2]{mrw}).

A \emph{locally compact quantum group} in the von Neumann algebraic setting
is a Hopf-von Neumann algebra $(M,\Delta)$ equipped with left and right
invariant weights.  As usual, we use $\Delta$ to turn $M_*$ into a Banach
algebra, and we write the product by juxtaposition.
We shall ``work on the left''; so using the left
invariant weight, we build the GNS space $H$, and a multiplicative
unitary $W$ acting on $L^2(\G)\otimes L^2(\G)$ (of course, the existence
of a right weight is needed to show that $W$ is unitary).  There is a
(in general unbounded) antipode $S$ which admits a ``polar decomposition''
$S=R \tau_{-i/2}$, where $R$ is the unitary antipode, and $(\tau_t)$ is
the scaling group.  There is a nonsingular positive operator $P$ which
implements $(\tau_t)$ as $\tau_t(x) = P^{it} x P^{-it}$.  Then $W$ is
\emph{manageable} with respect to $P$.

A manageable multiplicative unitary $W$ acting on $H\otimes H$
has, by definition, a nonsingular positive operator $P$, and an operator
$\tilde W$ acting on $H\otimes \overline{H}$ such that
\[ \big( W(\xi\otimes\alpha) \big| \eta\otimes\beta \big)
= \big( \tilde W(P^{-1/2}\xi\otimes\overline{\beta}) \big|
P^{1/2}\eta \otimes \overline{\alpha} \big), \]
for all $\alpha,\beta\in H$ and $\xi\in D(P^{-1/2}), \eta\in D(P^{1/2})$.
A word on notation: we work with left multiplicative unitaries, whereas
So{\l}tan and Woronowicz, in the conventions of \cite{kvvn}, work with
right multiplicative unitaries, and so we have translated everything to
the left.

Given such a $W$, the space $\{ (\iota\otimes\omega)W : \omega\in\mc B(H)_* \}$
is an algebra, and its closure is a C$^*$-algebra, say $A$.  There is a
coassociative map $\Delta:A\rightarrow M(A\otimes A)$ given by $\Delta(a)
= W^*(1\otimes a)W$.  If we formed $W$ from $(M,\Delta)$ with invariant weights,
then $A$ is $\sigma$-weakly dense in $M$, and the two definitions of $\Delta$
agree.  Similarly, $\{ (\omega\otimes\iota)W : \omega\in\mc B(H)_* \}$ is norm
dense in a C$^*$-algebra $\hat A$, and defining $\hat\Delta(\hat a) =
\hat W^*(1\otimes \hat a)\hat W$, we get a non-degenerate
$*$-homomorphism $\hat\Delta:\hat A \rightarrow M(\hat A\otimes \hat A)$,
where here $\hat W = \Sigma W^*\Sigma$.  If we started with $(M,\Delta)$ having
invariant weights, then we can construct invariant weights on $(\hat A'',
\hat\Delta)$.  The unitary $W$ is in the multiplier algebra $M(A\otimes \hat A)
\subseteq \mc B(H\otimes H)$.

When $W$ is a manageable multiplicative unitary, we can still form $S,R$
and $(\tau_t)$ with the usual properties.
The antipode $S$ has elements of the form $(\iota\otimes\omega)W$ as a core,
and $S((\iota\otimes\omega)W) = (\iota\otimes\omega)(W^*)$.  Then also $S
= R \tau_{-i/2}$ and $(\tau_t)$ is again implemented by $P$ (the same map
which appears in the definition of ``manageable'').

There is a more general notation of a \emph{modular multiplicative
unitary}, see \cite{sw}, which is more natural in certain examples.
However, at the cost of changing our space $H$, we can recover $(A,\Delta)$
from a different, but related, manageable multiplicative unitary.  Indeed,
in \cite{sw2}, it is shown that if $(A,\Delta)$ is given by some
modular multiplicative unitary, then $(\hat A,\hat\Delta)$, the $\sigma$-weak
topologies on $A$ and $\hat A$, the image of $W$ in $M(A\otimes \hat A)$,
and all the maps $S,R,(\tau_t),\hat S,\hat R$ and $(\hat\tau_t)$,
are independent of the particular choice of modular multiplicative unitary
giving $(A,\Delta)$.  For this reason, we shall henceforth work only with
manageable multiplicative unitaries (but of course all our results hold in
the modular case as well).

We write $\G$ for an abstract object to be thought of as a quantum group.
We write $C_0(\G), L^\infty(\G)$ and $L^1(\G)$ for $A,M$ and $M_*$
(and similar for the dual objects); as mentioned in the previous paragraph,
these are well-defined.  We also write $L^2(\G)$ for $H$, but be aware that if
$\G$ is given by a modular or manageable multiplicative unitary, then there is
some arbitrary choice involved in $L^2(\G)$.  If $\G$ has invariant weights,
then these weights unique up to a constant, and so $L^2(\G)$ is unique.

This concludes our brief summary; we shall develop further theory as and when
we need it.  We finish this section with one of our major tools-- that
``invariants are constant''.  Notice that, by using
the unitary antipode, we could replace $y_{13}$ by $y_{23}$ in the following;
but we shall have no need of this variant.

\begin{theorem}\label{thm:invcont}
For any $\G$ and a von Neumann algebra $N$, if $x,y\in L^\infty(\G)\vnten N$
satisfy $(\Delta\otimes\iota)x = y_{13}$, then $x=y\in \mathbb C\vnten N$.
\end{theorem}
\begin{proof}
We shall prove this when $N=\mathbb C$, the general case comes from
considering $(\iota\otimes\omega)x$ and $(\iota\otimes\omega)y$, as
$\omega\in N_*$ varies.  For locally compact quantum groups, this was
shown in \cite[Lemma~4.6]{ari}, compare also \cite[Result~5.13]{kv}.
For general $\G$, \cite[Theorem~2.6]{mrw} shows that if $a,b\in\mc B(L^2(\G))$
with $W^*(1\otimes a)W = b\otimes 1$ (working with left multiplicative
unitaries) then $a=b\in\mathbb C1$, and this immediately implies the result.
\end{proof}

The following result is, when we have invariant weights, well-known
to experts, but proofs can be hard to find (for example, \cite{jnr} cites
\cite[Proposition~3.5]{vv2}, which uses Tomita-Takesaki theory;
compare also the proof of \cite[Proposition~5.13]{vv}).  The
argument here was suggested to us by the anonymous referee of an earlier
version of this paper.

\begin{theorem}\label{thm:gen_all_bh}
Let $\G$ be a quantum group.  Then the linear span of $L^\infty(\hat\G)
L^\infty(\G)' = \{ \hat x y' : \hat x\in L^\infty(\hat\G),
y'\in L^\infty(\G)' \}$ is $\sigma$-weakly dense in $\mc B(L^2(\G))$.
\end{theorem}
\begin{proof}
Let $N$ be the $\sigma$-weak closed linear span of $L^\infty(\hat\G)
L^\infty(\G)'$, and let $P$ be the von Neumann algebra generated by $N$.
We first claim that $N=P$; we need only show that $P\subseteq N$.  In fact,
it suffices to show that $L^\infty(\G)' L^\infty(\hat\G) \subseteq N$,
as then $N$ will be a bimodule over both $L^\infty(\hat\G)$ and
$L^\infty(\G)'$, and hence will contain the algebra generated by
$L^\infty(\hat\G)$ and $L^\infty(\G)'$.

Let $y'\in L^\infty(\G)'$, let $\xi,\eta\in L^2(\G)$ and set
$x=(\iota\otimes\omega_{\xi,\eta})(\hat W) \in L^\infty(\hat\G)$.
Letting $(e_i)$ be an orthonormal basis for $L^2(\G)$, we have that
\begin{align*} y'x &= y' (\iota\otimes\omega_{\xi,\eta})(\hat W)
= (\iota\otimes\omega_{\xi,\eta})(\hat W \hat W^*(y'\otimes 1)\hat W) \\
&= \sum_i (\iota\otimes\omega_{e_i,\eta})(\hat W)
(\iota\otimes\omega_{\xi,e_i})(\hat W^*(y'\otimes 1)\hat W), \end{align*}
the sum converging $\sigma$-weakly.  However, for $\omega\in \mc B(L^2(\G))_*$
and $y\in L^\infty(\G)$,
\begin{align*} (\iota\otimes\omega)(\hat W^*(y'\otimes 1)\hat W) y
&= (\omega\otimes\iota)(W(1\otimes y')W^*(1\otimes y))
= (\omega\otimes\iota)(W(1\otimes y')\Delta(y) W) \\
&= (\omega\otimes\iota)(W\Delta(y)(1\otimes y') W)
= (\omega\otimes\iota)((1\otimes y)W(1\otimes y') W) \\
&= y (\iota\otimes\omega)(\hat W^*(y'\otimes 1)\hat W).
\end{align*}
It follows that $(\iota\otimes\omega)(\hat W^*(y'\otimes 1)\hat W) \in
L^\infty(\G)'$, and so
$y' (\iota\otimes\omega_{\xi,\eta})(\hat W) \in N$.  As such $x$ are
$\sigma$-weakly dense in $L^\infty(\hat\G)$, the claim follows.

Now let $z\in P' = L^\infty(\hat\G)' \cap L^\infty(\G)$.
Then $\Delta(z) = W^*(1\otimes z)W = 1\otimes z$ as $W\in L^\infty(\G)
\vnten L^\infty(\hat\G)$.  By applying Theorem~\ref{thm:invcont} to
$R(z)$ it follows that $z\in\mathbb C1$.  Thus $N=P=(\mathbb C1)''
= \mc B(L^2(\G))$ as required.
\end{proof}

\section{Multipliers of quantum groups}\label{sec:mult}

In this section, we review some of the ideas used by Junge, Neufang and
Ruan in \cite{jnr}.  We shall actually need some constructions coming from the
proofs in \cite{jnr} (and not just the statements of the results).  Rather than
just give sketch proofs, we instead give quick, full proofs, and take the
opportunity to show that some of their results also hold for quantum groups
coming from manageable multiplicative unitaries.  Further details and related
ideas can be found in \cite{daws1,daws2,hnr2,jnr}.

\begin{definition}
A completely bounded left multiplier of $L^1(\G)$ is a completely bounded map
$L_*:L^1(\G) \rightarrow L^1(\G)$ with
$L_*(\omega_1\omega_2) = L_*(\omega_1) \omega_2$ for
$\omega_1,\omega_2 \in L^1(\G)$.
\end{definition}

Such maps are also often called ``centralisers'' in the literature (and in
particular, in \cite{jnr}).  A simple calculation shows that a completely
bounded map $L_*:L^1(\G) \rightarrow L^1(\G)$ is a left multiplier if and
only if its adjoint $L = (L_*)^*$ satisfies $(L\otimes\iota)\Delta =
\Delta L$.

Let us make a few remarks about normal completely bounded maps.
As explained, for example, in the proof of \cite[Theorem~2.5]{HM}, as
$L$ is normal, we can find a normal $*$-representation $\pi:L^\infty(\G)
\rightarrow \mc B(H)$ for some Hilbert space $H$, and bounded maps
$P,Q:L^2(\G)\rightarrow H$, with $L(x) = P^*\pi(x)Q$ for each
$x\in L^\infty(\G)$.  By the structure theory for normal $*$-representations
(see \cite[Theorem~5.5, Chapter~IV]{tak1}) by adjusting $P$ and $Q$,
and may suppose that $H=L^2(\G)\otimes H'$ for some Hilbert space $H'$,
and that $\pi(x) = x\otimes 1$.  For example, it then follows that
\[ (L\otimes\iota)(\hat W) = (P^*\otimes 1)\hat W_{13}(Q\otimes 1). \]
As $\hat W\in M(\mc B_0(L^2(\G)) \otimes C_0(\G))$, it follows easily from
this that also $(L\otimes\iota)(\hat W)\in M(\mc B_0(L^2(\G)) \otimes C_0(\G))$,
a fact we shall use in the following proof.

The following is a short unification of (the left version of)
\cite[Corollary~4.4]{jnr} (compare \cite[Theorem~4.10]{jnr})
and \cite[Theorem~4.2]{daws1}; we make use
of the ``invariants are constant'' technique.

Remember that the \emph{left regular representation} of $\G$ is the
injective homomorphism $\lambda:L^1(\G)\rightarrow C_0(\hat\G);
\omega \mapsto (\omega\otimes\iota)(W)$.  By duality we define
$\hat\lambda:L^1(\hat\G)\rightarrow C_0(\G); \hat\omega \mapsto
(\hat\omega\otimes\iota)(\hat W) = (\iota\otimes\hat\omega)(W^*)$.

\begin{proposition}\label{prop:mults_rep}
Let $L_*$ be a completely bounded left multiplier of $L^1(\hat\G)$.
There is $a\in M(C_0(\G))$ with $a \hat\lambda(\hat\omega)
= \hat\lambda(L_*(\hat\omega))$ for $\hat\omega\in L^1(\G)$, or equivalently,
with $(1\otimes a)\hat W = (L\otimes\iota)(\hat W)$.
\end{proposition}
\begin{proof}
That $a \hat\lambda(\hat\omega) = \hat\lambda(L_*(\hat\omega))$ for each
$\hat\omega\in L^1(\G)$, if and only if $(1\otimes a)\hat W = 
(L\otimes\iota)(\hat W)$ follows easily from the definition of $\hat\lambda$.
Consider now
\begin{align*}
(\hat\Delta\otimes\iota)\big( (L\otimes\iota)(\hat W) \hat W^* \big)
&= \big( (L\otimes\iota\otimes\iota) (\hat\Delta\otimes\iota)(\hat W) \big)
(\hat\Delta\otimes\iota)(\hat W^*) \\
&= \big( (L\otimes\iota\otimes\iota) (\hat W_{13} \hat W_{23}) \big)
\hat W_{23}^* \hat W_{13}^* \\
&= (L\otimes\iota)(\hat W)_{13} \hat W_{13}^*, \end{align*}
where we have used that $\hat\Delta L = (L\otimes\iota)\hat\Delta$,
that $\hat\Delta$ is a $*$-homomorphism, and that
$(\hat\Delta\otimes\iota)(\hat W) = \hat W_{13} \hat W_{23}$.
By Theorem~\ref{thm:invcont}, it follows that there is $a\in L^\infty(\G)$
with $(L\otimes\iota)(\hat W) \hat W^* = 1\otimes a$.  However, as
$\hat W \in M(\mc B_0(L^2(\G)) \otimes C_0(\G))$, we see that
\[ 1 \otimes a = (L\otimes\iota)(\hat W) \hat W^* \in M(\mc B_0(L^2(\G))
\otimes C_0(\G)), \]
from which it follows immediately that $a\in M(C_0(\G))$.
\end{proof}

In the language of \cite{daws1}, the previous lemma says that $L_*$ is
``represented''; in the language of \cite{jnr}, the element $a$ is the
``multiplier'' associated to the ``centraliser'' $L_*$.

\begin{proposition}\label{prop:one}
Let $L_*$ be a completely bounded left multiplier of $L^1(\hat\G)$.
There is a completely bounded, normal map $\Phi : \mc B(L^2(\G)) \rightarrow
\mc B(L^2(\G))$ which extends $L$, and which is a $L^\infty(\G)'$-bimodule
map.  Indeed, $\Phi$ satisfies
\[ 1\otimes \Phi(x) = \hat W \big( (L\otimes\iota)(\hat W^*(1\otimes x)\hat W)
\big) \hat W^* \qquad (x\in\mc B(L^2(\G))). \]
\end{proposition}
\begin{proof}
We closely follow \cite[Proposition~4.3]{jnr}, while translating ``to
the left'' and using that ``invariants are constant''.  Define
\[ T:\mc B(L^2(\G)) \rightarrow L^\infty(\hat\G) \vnten \mc B(L^2(\G)),
\quad T(x) = \hat W \big( (L\otimes\iota)(\hat W^*(1\otimes x)\hat W) \big)
\hat W^*. \]
We now perform a similar calculation to that in the previous lemma:
\begin{align*}
(\hat\Delta\otimes\iota)T(x) &=
\hat W_{13} \hat W_{23} (L\otimes\iota\otimes\iota)\big(
(\hat\Delta\otimes\iota)(\hat W^*(1\otimes x)\hat W)
\big) \hat W_{23}^* \hat W_{13}^* \\
&= \hat W_{13} \hat W_{23} (L\otimes\iota\otimes\iota)\big(
\hat W_{23}^* \hat W_{13}^* (1\otimes 1\otimes x) \hat W_{13} \hat W_{23}
\big) \hat W_{23}^* \hat W_{13}^* \\
&= \hat W_{13} (L\otimes\iota\otimes\iota)
\big( \hat W_{13}^* (1\otimes 1\otimes x) \hat W_{13} \big) \hat W_{13}^* \\
&= T(x)_{13}.
\end{align*}
So by Theorem~\ref{thm:invcont}, there is $\Phi(x)\in\mc B(L^2(\G))$ with
$T(x) = 1\otimes \Phi(x)$.  It is easy to see that $\Phi$ is completely
bounded and normal.

For $x\in L^\infty(\hat G)$
\[ 1\otimes\Phi(x) = T(x) = \hat W( (L\otimes \iota)\hat\Delta(x) )\hat W^*
= \hat W\hat\Delta(L(x))\hat W^* = 1\otimes L(x), \]
and so $\Phi$ extends $L$.  For $y,z\in L^\infty(\G)'$, as $\hat W
\in L^\infty(\hat\G) \vnten L^\infty(\G)$, it is easy to see that
\[ T(yxz) = (1\otimes y)T(x)(1\otimes z) \qquad (x\in\mc B(L^2(\G))), \]
and so $\Phi(yxz) = y\Phi(x)z$ as required.
\end{proof}

In the language of \cite{jnr}, we have thus constructed a map
from the set of completely bounded left multipliers of $L^1(\hat\G)$ to
$\mc{CB}^{\sigma, L^\infty(\hat\G)}_{L^\infty(\G)'}
(\mc B(L^2(\G)))$.  It seems that, to continue with the arguments of
\cite{jnr}, we start to need to use arguments that involve the relative
position of $L^\infty(\G)$ and its commutant in $\mc B(L^2(\G))$.
In particular, to show
that every $\Phi \in \mc{CB}^{\sigma, L^\infty(\hat\G)}_{L^\infty(\G)'}
(\mc B(L^2(\G)))$ comes from a left multiplier would require us to know that
$L^\infty(\G) \cap L^\infty(\hat\G) = \mathbb C$ (at least if one is
following the proof of \cite[Proposition~3.2]{jnr}), and we have no proof of
this in the Manageable Multiplicative Unitary setting.

\section{Multipliers coming from invertible corepresentations}
\label{sec:corep-mult}

In this section, we show, rather explicitly, how corepresentations and
``universal'' quantum groups give rise to completely bounded and completely
positive multipliers.  

A \emph{corepresentation} of $\G$ shall be, for us, an element
$U\in L^\infty(\G) \vnten \mc B(H)$ with $(\Delta\otimes\iota)(U)
= U_{13} U_{23}$ (so, we don't assume that $U$ is unitary).
We state the following in a little generality, but note that it obviously
applies to unitary corepresentations.  Similar ideas are explored in
\cite[Section~6]{daws1}, and for Kac algebras, with much less emphasis on
corepresentation theory, see \cite{kr}.

\begin{proposition}\label{prop:inv_corep_mult}
Let $U$ be a corepresentation of $\G$, and suppose there is
$V\in\mc B(L^2(\G)\otimes H)$ with $VU^* = 1$ (that is, $U$ has a right
inverse).  For each $\alpha,\beta\in H$, there is a completely bounded
left multiplier of $L^1(\hat\G)$ represented by
$a = (\iota\otimes\omega_{\alpha,\beta})(U^*)$.
If $U^*$ is an isometry (so we may take $V = U$) and $\alpha=\beta$, then
the multiplier is completely positive.
\end{proposition}
\begin{proof}
We have that $(\Delta\otimes\iota)(U^*) = U^*_{23} U^*_{13}$, or equivalently,
$W^*_{12} U^*_{23} = U^*_{23} U^*_{13} W_{12}^*$.  Thus also
$V_{23} W^*_{12} U^*_{23} = U^*_{13} W_{12}^*$, and using that $\hat W=
\Sigma W^* \Sigma$, it follows that $V_{13} \hat W_{12} U^*_{13} =
U^*_{23} \hat W_{12}$.  Thus define $L:L^\infty(\hat\G) \rightarrow
\mc B(L^2(\G))$ by
\[ L(\hat x) = (\iota\otimes\omega_{\alpha,\beta})(V(\hat x\otimes 1)U^*)
\qquad (\hat x\in L^\infty(\hat\G)). \]
Clearly $L$ is a normal, completely bounded map.
Then immediately we see that $(L\otimes\iota)(\hat W) = (1\otimes a)\hat W$, 
and it is now easy to see (compare \cite[Proposition~2.3]{daws1}) that $L$
maps into $L^\infty(\hat\G)$, and that $L$ is the adjoint of a left
multiplier on $L^1(\hat\G)$, represented by $a$.

When $U^*$ is an isometry, $V=U$ and $\alpha=\beta$, clearly $L$ is
completely positive.
\end{proof}

For $U,V$ as in the proposition, we could weaken the condition on $U$
to asking that $U\in\mc B(L^2(\G)\otimes H)$ with $W_{12}^* U_{23} W_{12}
= U_{13} U_{23}$.  Then, arguing as in \cite[Page~142]{woro}, we see that
$U_{13} = W_{12}^* U_{23} W_{12} V_{23}^* \in M(C_0(\G) \otimes \mc B_0(L^2(\G))
\otimes \mc B_0(H))$, and so $U\in M(C_0(\G)\otimes\mc B_0(H))$, in particular,
$U$ is a corepresentation in our sense.

Let us just remark that if also $V$ is a corepresentation, then consider
forming $\Phi$ as in Section~\ref{sec:mult}, using the $L$ given as in
the proposition.  So, for $x\in\mc B(L^2(\G))$,
\begin{align*} 1\otimes\Phi(x) &=
(\iota\otimes\omega_{\alpha,\beta}\otimes\iota)
   \big( \hat W_{13} V_{12} \hat W^*_{13} (1\otimes 1 \otimes x) \hat W_{13}
   U_{12}^* \hat W_{13}^* \big) \\
&= (\iota\otimes\iota\otimes\omega_{\alpha,\beta}) \big(
\hat W_{12} V_{13} \hat W^*_{12} (1\otimes x \otimes 1) \hat W_{12}
   U_{13}^* \hat W_{12}^* \big).
\end{align*}
Now, $\hat W(a\otimes 1)\hat W^* = \Sigma W^*(1\otimes a)W \Sigma
= \Sigma\Delta(a)\Sigma$ for $a\in L^\infty(\G)$, and so
\[ 1\otimes\Phi(x) = (\iota\otimes\iota\otimes\omega_{\alpha,\beta}) \big(
V_{23} V_{13} (1\otimes x\otimes 1) U^*_{13} U^*_{23} \big)
= 1 \otimes (\iota\otimes\omega_{\alpha,\beta})
\big( V (x\otimes 1) U^* \big). \]
Hence $\Phi$ has the same ``defining formula'' as $L$.

\subsection{Links with universal quantum groups}\label{sec:universal_link}

Universal quantum groups are constructed in \cite{kus} and \cite[Section~5]{sw}.
We write $C_0^u(\hat\G)$ for the universal dual of $C_0(\G)$.  For us,
the important properties are that:
\begin{itemize}
\item There is a coassociative non-degenerate $*$-homomorphism
$\hat\Delta_u:C_0^u(\hat\G) \rightarrow M(C_0^u(\hat\G)\otimes C_0^u(\hat\G))$;
\item There is a surjective $*$-homomorphism $\hat\pi_u:C_0^u(\hat\G)
\rightarrow C_0(\hat\G)$ with $\hat\Delta\hat\pi_u =
(\hat\pi_u\otimes\hat\pi_u) \hat\Delta_u$;
\item There is a unitary corepresentation $\mc W\in M(C_0(\G)\otimes
C_0^u(\hat\G))$ of $C_0(\G)$ such that $(\iota\otimes\hat\pi_u)\mc W=W$
and $(\iota\otimes\hat\Delta_u)\mc W = \mc W_{13} \mc W_{12}$.
\item The space $\{(\omega\otimes\iota)(\mc W) : \omega\in L^1(\G)\}$
is dense in $C_0^u(\hat\G)$.
\item There is a bijection between unitary corepresentations $U$ of
$C_0(\G)$ and non-degenerate $*$-homomorphisms $\pi:C_0^u(\hat\G)
\rightarrow\mc B(H)$ given by the relation that $U = (\iota\otimes\pi)(\mc W)$.
\end{itemize}
Note that our $\mc W$ is denoted by $\hat{\mc V}$ in the notation of \cite{kus};
and is the ``left analogue'' of $\mathbb W$ in the notation of \cite{sw}.

The map $\hat\pi_u^*:L^1(\hat\G)\rightarrow C_0^u(\hat\G)^*$ is an isometry and
an algebra homomorphism.  We know (see for example \cite[Proposition~8.3]{daws2})
that this identifies $L^1(\hat\G)$ with an ideal in $C_0^u(\hat\G)^*$, and
hence that members of $C_0^u(\hat\G)^*$ induce multipliers on $L^1(\hat\G)$.
Let us make links with Proposition~\ref{prop:inv_corep_mult}.

\begin{proposition}\label{prop:mults_from_uni_dual}
Let $U$ be a unitary corepresentation of $\G$ on $H$, and let
$\alpha,\beta\in H$.  Let $\pi$ be the $*$-representation of $C_0^u(\hat\G)$
on $H$ associated with $U$.  Then the multiplier represented by
$(\iota\otimes\omega_{\alpha,\beta})(U^*)$ is given by left multiplication
by $\mu = \omega_{\alpha,\beta}\circ\pi\in C_0^u(\hat\G)^*$.
\end{proposition}
\begin{proof}
Let $L:L^\infty(\hat\G) \rightarrow L^\infty(\hat\G)$ be the adjoint of
the completely bounded left multiplier represented by
$a = (\iota\otimes\omega_{\alpha,\beta})(U^*)$, as constructed in
Proposition~\ref{prop:inv_corep_mult}.  Then $(L\otimes\iota)(\hat W)
= (1\otimes a)\hat W$, or equivalently, $(\iota\otimes L)(W^*)
= (a\otimes 1)W^*$.

Define $L^\dagger(\hat x) = L(\hat x^*)^*$ for $\hat x\in L^\infty(\hat\G)$,
so $L^\dagger$ is a normal, completely bounded map on $L^\infty(\hat\G)$.
For any von Neumann algebra $M$ and $X\in M\vnten L^\infty(\hat\G)$,
we see that $(\iota\otimes L^\dagger)(X^*) = (\iota\otimes L)(X)^*$.
In particular, it follows that $(L^\dagger\otimes\iota)\hat\Delta
= \hat\Delta L^\dagger$, and so $L^\dagger$ is the adjoint of a completely
bounded left multiplier of $L^1(\hat\G)$, represented by $b$ say.
The proof of Proposition~\ref{prop:inv_corep_mult} shows that
$b = (\iota\otimes\omega_{\beta,\alpha})(U^*)$.

Given $\hat\omega\in L^1(\hat\G)$, we wish to show that $\mu\hat\pi_u^*(\hat\omega)
= \hat\pi_u^*(L_*(\hat\omega))$.  Let $\omega\in L^1(\G)$ and set
$x = (\omega\otimes\iota)\mc W \in C_0^u(\hat\G)$.  Then
\begin{align*}
\ip{\mu\hat\pi_u^*(\hat\omega)}{x}
&= \ip{\mu\otimes\hat\pi_u^*(\hat\omega)}{\hat\Delta_u((\omega\otimes\iota)\mc W)}
= \ip{\omega\otimes\mu\otimes\hat\pi_u^*(\hat\omega)}{\mc W_{13} \mc W_{12}} \\
&= \ip{\omega\otimes\omega_{\alpha,\beta}\otimes\hat\omega}
   {(\iota\otimes\hat\pi_u)(\mc W)_{13} (\iota\otimes\pi)(\mc W)_{12}}
= \ip{W_{13} U_{12}}{\omega\otimes\omega_{\alpha,\beta}\otimes\hat\omega},
\end{align*}
and also
\[ \ip{\hat\pi_u^*(L_*(\hat\omega))}{x}
= \ip{(\iota\otimes\hat\pi_u)\mc W}{\omega\otimes L_*(\hat\omega)}
= \ip{W}{\omega\otimes L_*(\hat\omega)}
= \ip{(\iota\otimes L)(W)}{\omega\otimes\hat\omega}. \]
Now, $(\iota\otimes L)(W) = (\iota\otimes L^\dagger)(W^*)^*
= ((b\otimes 1)W^*)^* = W(b^*\otimes 1)$, and so, using that
$b^* = (\iota\otimes\omega_{\alpha,\beta})(U)$, we have that
\[ \ip{\hat\pi_u^*(L_*(\hat\omega))}{x}
= \ip{W(b^*\otimes 1)}{\omega\otimes\hat\omega}
= \ip{W_{13} U_{12}}{\omega\otimes \omega_{\alpha,\beta} \otimes\hat\omega}. \]
As such $x$ are dense in $C_0^u(\hat\G)$, the proof is complete.
\end{proof}

In particular, taking $U=\mc W$, we see that every positive functional on
$C_0^u(\hat\G)$ induces a completely positive left multiplier of $L^1(\hat\G)$.
The main result of this paper is to show that the converse is also true.

\section{Completely positive multipliers}\label{sec:cpmult}

In this section, we study completely \emph{positive} multipliers of
$L^1(\hat\G)$.  Motivated by Proposition~\ref{prop:one}, we will first
study completely positive normal maps on $\mc B(L^2(\G))$.  As
$\mc B_0(L^2(\G))$ is $\sigma$-weakly-dense in $\mc B(L^2(\G))$, it
suffices to consider completely positive maps $\mc B_0(L^2(\G))
\rightarrow \mc B(L^2(\G))$.  The following is simply the Stinespring
construction, tailored to this specific situation.  We give the details,
as they are central to our argument.  For $\xi,\eta\in
L^2(\G)$, let $\theta_{\xi,\eta}\in\mc B_0(L^2(\G))$ be the rank-one
operator $\alpha \mapsto (\alpha|\eta) \xi$.

Let $\Phi:\mc B_0(L^2(\G)) \rightarrow \mc B(L^2(\G))$ be a completely
positive map.  We remark that $\Phi$ has a unique normal extension to
a normal completely positive map $\mc B(L^2(\G)) \rightarrow \mc B(L^2(\G))$,
which we shall also denote by $\Phi$.
Let $H$ be the completion of the algebraic tensor product
$\overline{L^2(\G)} \odot L^2(\G)$ for the pre-inner-product
\[ \big( \overline{\xi}\otimes\alpha \big| \overline{\eta}\otimes\beta \big)_H
= \big( \Phi(\theta_{\eta,\xi}) \alpha | \beta \big). \]
That this is a \emph{positive} sesquilinear form follows from the fact
that $\Phi$ is completely positive (compare with \cite[Theorem~3.6]{tak1}
for example).  We shall write $\nu_{\alpha,\xi}$ for the equivalence
class of $\overline{\xi}\otimes\alpha$ in $H$; see the end of the following
paragraph for an explanation of this notation.

Let $(e_i)$ be an orthonormal basis of $L^2(\G)$, and define
\[ V:L^2(\G) \rightarrow L^2(\G)\otimes H;
\quad \alpha \mapsto \sum_i e_i \otimes \nu_{\alpha,e_i}. \]
This makes sense, that is, the sum converges, as
\[ \sum_i \|\nu_{\alpha,e_i}\|^2_H
= \sum_i (\Phi(\theta_{e_i,e_i})\alpha|\alpha)
= (\Phi(1)\alpha|\alpha). \]
Then, for $\alpha,\beta,\xi,\eta\in L^2(\G)$,
\begin{align*}
(V^*(\theta_{\xi,\eta}\otimes 1)V\alpha|\beta)
&= \sum_{i,j} (\theta_{\xi,\eta}(e_i) \otimes \nu_{\alpha,e_i} |
e_j \otimes \nu_{\beta,e_j} )_H \\
&= \sum_i (e_i|\eta) (\xi|e_j) (\Phi(\theta_{e_j,e_i})\alpha|\beta)
= (\Phi(\theta_{\xi,\eta})\alpha|\beta).
\end{align*}
Thus we have a Stinespring dilation of $\Phi$.  Now let $(f_i)$ be
an orthonormal basis of $H$, and define a family $(a_i)$ in $\mc B(L^2(\G))$
by setting
\[ V(\alpha) = \sum_i a_i(\alpha) \otimes f_i \in L^2(\G)\otimes H
\qquad (\alpha\in L^2(\G)). \]
It follows that $\Phi(x) = \sum_i a_i^* x a_i$ for each $x\in\mc B(L^2(\G))$.
Furthermore, for $\xi,\eta,\alpha,\beta\in L^2(\G)$,
\[ (\nu_{\alpha,\xi}|\nu_{\beta,\eta})_H
= \sum_i (a_i^*\theta_{\eta,\xi}a_i\alpha|\beta)
= \sum_i (a_i\alpha|\xi) (\eta|a_i\beta)
= \sum_i \ip{a_i}{\omega_{\alpha,\xi}}
\overline{ \ip{a_i}{\omega_{\beta,\eta}} }, \]
and so $\nu_{\alpha,\xi} = \sum_i \ip{a_i}{\omega_{\alpha,\xi}} f_i$ in $H$.
This explains the choice of notation $\nu_{\alpha,\xi}$, which is deliberately
reminiscent of $\omega_{\alpha,\xi}$.

\begin{proposition}\label{prop:two}
In the above setting, suppose further that $M$ is a von Neumann algebra
on $L^2(\G)$, and that $\Phi$ is an $M$-bimodule map.  Then
$(x\otimes 1)V = Vx$ for each $x\in M$, and $a_i\in M'$ for each $i$.
\end{proposition}
\begin{proof}
For readability, we drop the $\nu$ notation in this proof.
For $x\in M$ and $\xi,\eta,\alpha,\beta\in L^2(\G)$,
\[ ( \overline{x^*(\xi)}\otimes\alpha|\overline{\eta}\otimes\beta)_H
= (\Phi(\theta_{\eta,\xi}x)\alpha|\beta)
= (\Phi(\theta_{\eta,\xi})x\alpha|\beta)
= (\overline\xi\otimes x(\alpha)|\overline\eta\otimes\beta)_H. \]
Thus $\overline{x^*(\xi)}\otimes\alpha = \overline\xi\otimes x(\alpha)$ in $H$.
It follows that
\begin{align*} V(x(\alpha)) &= \sum_i e_i \otimes
   (\overline{e_i}\otimes x(\alpha))
= \sum_i e_i \otimes (\overline{x^*(e_i)}\otimes \alpha)
= \sum_{i,j} e_i \otimes ( \overline{(x^*(e_i)|e_j)e_j} \otimes \alpha) \\
&= \sum_{i,j} (x(e_j)|e_i) e_i \otimes ( \overline{e_j} \otimes \alpha)
= \sum_{j} x(e_j) \otimes ( \overline{e_j} \otimes \alpha)
= (x\otimes 1)V(\alpha), \end{align*}
remembering that $H$ is the completion of $\overline{L^2(\G)}\otimes L^2(\G)$.
It now follows that $xa_i = a_ix$ for each $i$, and so as $x\in M$ was
arbitrary, $a_i \in M'$ for each $i$.
\end{proof}

The previous result (in the more general completely bounded setting) is
well-known, see for example \cite[Theorem~3.1]{smith} and unpublished work of
Haagerup.  However, the actual construction will be central to our arguments.

\subsection{Constructing a corepresentation}

For the remainder of this section, fix a completely positive left multiplier
on $L^1(\hat\G)$.  Form $\Phi:\mc B(L^2(\G)) \rightarrow \mc B(L^2(\G))$
using Proposition~\ref{prop:one} applied to this multiplier, and apply the
construction of the previous section to find $H$ and
$V:L^2(\G) \rightarrow L^2(\G)\otimes H$.
Fixing an orthonormal basis $(f_i)$ for $H$, we find $a_i$ such that
$\Phi(x) = \sum_i a_i^* x a_i$ for each $x\in\mc B(L^2(\G))$.
By Proposition~\ref{prop:two}, we see that $(x\otimes 1)V = Vx$ for
each $x\in L^\infty(\G)'$, equivalently, that $a_i\in L^\infty(\G)$
for each $i$.

\begin{proposition}
There is a unique isometry $U^*$ on $L^2(\G)\otimes H$ which
satisfies
\[ U^*\big( \xi \otimes \nu_{\alpha,\eta} \big)
= \sum_i (\omega_{\alpha,\eta}\otimes\iota)\Delta(a_i)\xi\otimes f_i, \]
for all $\xi,\eta,\alpha\in L^2(\G)$.
\end{proposition}
\begin{proof}
We know that $\Phi(x) = \sum_i a_i^* x a_i$ for
$x\in L^\infty(\hat\G)$.  We now use Proposition~\ref{prop:one},
which tells us that
\[ 1\otimes\Phi(x) = \sum_i \hat W (a_i^*\otimes 1)\hat W^*(1\otimes x)\hat W
(a_i\otimes 1) \hat W^* \qquad (x\in\mc B_0(L^2(\G))). \]
For $\xi_1,\eta_1,\alpha_1,\xi_2,\eta_2,\alpha_2\in L^2(\G)$, we hence
have that
\begin{align*}
\big( \xi_1\otimes \nu_{\alpha_1,\eta_1} \big| &
   \xi_2\otimes \nu_{\alpha_2,\eta_2} \big)_{L^2(\G)\otimes H}
= \big( (1\otimes\Phi(\theta_{\eta_2,\eta_1})) (\xi_1\otimes\alpha_1)
\big| \xi_2\otimes\alpha_2 \big) \\
&= \sum_i \big( \hat W (a_i^*\otimes 1)\hat W^*(1\otimes \theta_{\eta_2,\eta_1})
\hat W (a_i\otimes 1) \hat W^* (\xi_1\otimes\alpha_1)
\big| \xi_2\otimes\alpha_2 \big) \\
&= \sum_i \big( (1\otimes \theta_{\eta_2,\eta_1}) \Sigma \Delta(a_i) \Sigma
(\xi_1\otimes\alpha_1) \big|
\Sigma \Delta(a_i) \Sigma (\xi_2\otimes\alpha_2) \big) \\
&= \sum_i \big( (\theta_{\eta_2,\eta_1}\otimes 1) \Delta(a_i)
(\alpha_1\otimes\xi_1) \big| \Delta(a_i) (\alpha_2\otimes\xi_2) \big) \\
&= \sum_i \big( (\omega_{\alpha_1,\eta_1}\otimes\iota)\Delta(a_i) \xi_1 \big|
(\omega_{\alpha_2,\eta_2}\otimes\iota)\Delta(a_i) \xi_2 \big),
\end{align*}
using that $\hat W(a\otimes 1)\hat W^* = \Sigma \Delta(a) \Sigma$
for $a\in L^\infty(\G)$.  It follows immediately that $U^*$ exists and
is an isometry; uniqueness follows as vectors of the form
$\xi\otimes\nu_{\alpha,\eta}$ are linearly dense in $L^2(\G)\otimes H$.
\end{proof}

As $\nu_{\alpha,\eta} = \sum_i \ip{a_i}{\omega_{\alpha,\eta}} f_i$,
by using linearity and continuity, we see that also
\[ U^*\Big( \xi\otimes \sum_i \ip{a_i}{\omega} f_i \Big)
= \sum_i (\omega\otimes\iota)\Delta(a_i) \xi \otimes f_i
\qquad (\xi\in L^2(\G), \omega\in L^1(\G)). \]

\begin{proposition}\label{prop:iscorep}
The operator $U$ is a member of $L^\infty(\G) \vnten \mc B(H)$, and
is a corepresentation, that is, $(\Delta\otimes\iota)U = U_{13} U_{23}$.
\end{proposition}
\begin{proof}
Let $x\in L^\infty(\G)'$, so for $\xi,\alpha,\beta\in L^2(\G)$,
\begin{align*} U^*(x\xi\otimes \nu_{\alpha,\beta})
&= \sum_i (\omega_{\alpha,\beta}\otimes\iota)\Delta(a_i)x\xi\otimes f_i
= \sum_i x (\omega_{\alpha,\beta}\otimes\iota)\Delta(a_i)\xi\otimes f_i \\
&= (x\otimes 1)U^*(\xi\otimes \nu_{\alpha,\beta}). \end{align*}
Thus $U^* \in (L^\infty(\G)'\vnten \mathbb C)' =
L^\infty(\G) \vnten \mc B(H)$, and of course the same is true of $U$.

We shall prove that $(\Delta\otimes\iota)(U^*) = U_{23}^* U_{13}^*$.  It is
easy to see that this is equivalent to $\pi:L^1(\G)\rightarrow\mc B(H);
\omega\mapsto (\omega\otimes\iota)(U^*)$ being an anti-homomorphism of the
Banach algebra $L^1(\G)$.  However, notice that for $\omega_1,\omega_2\in L^1(\G)$
and $\xi,\eta\in L^2(\G)$,
\begin{align*}
\Big( \pi(\omega_{\xi,\eta}) \sum_i \ip{a_i}{\omega_1}f_i \Big|
   \sum_j \ip{a_j}{\omega_2}f_j \Big)
&= \Big( U^*\Big( \xi\otimes\sum_i \ip{a_i}{\omega_1}f_i\Big) \Big|
   \eta\otimes\sum_j \ip{a_j}{\omega_2}f_j \Big) \\
&= \Big( \sum_i (\omega_1\otimes\iota)\Delta(a_i)\xi \otimes f_i \Big|
   \eta\otimes\sum_j \ip{a_j}{\omega_2}f_j \Big) \\
&= \Big( \sum_i \ip{a_i}{\omega_1 \omega_{\xi,\eta}} f_i \Big|
   \sum_j \ip{a_j}{\omega_2} f_j \Big).
\end{align*}
Thus
\[ \pi(\omega)\Big( \sum_i \ip{a_i}{\omega'}f_i \Big)
= \sum_i \ip{a_i}{\omega'\omega}f_i
\qquad (\omega,\omega'\in L^1(\G)), \]
and it is now immediate that $\pi$ is an anti-homomorphism.
\end{proof}

We remark that we can view $H$ as being a completion of $L^1(\G)$, where
we identify $\omega\in L^1(\G)$ with $\sum_i \ip{a_i}{\omega} f_i\in H$.
Then $\pi$ in the above proof (that is, the anti-homomorphism from $L^1(\G)$
to $\mc B(H)$ induced by $U^*$) is simply the map $\pi(\omega):\omega'
\mapsto\omega'\omega$.

We now finish the argument by showing that $U$ is unitary.

\begin{lemma}\label{lem:image_ustar}
The closed image of $U^*$ is equal to the closed linear span of
$\{ (\hat a\otimes 1)V(\xi) : \xi\in L^2(\G), \hat a\in C_0(\hat\G) \}$.
In particular, the image of $U^*$ contains the image of $V$, and so $U^*UV=V$.
\end{lemma}
\begin{proof}
Let $\xi_1,\xi_2,\eta\in L^2(\G)$, and let $\sum_i \xi_i \otimes f_i \in
L^2(\G) \otimes H$, and observe that
\begin{align*} \Big( U^*\big( \xi_1 \otimes \sum_i (a_i\xi_2|\eta) f_i \big)
   \Big| \sum_j \xi_j \otimes f_j \Big)
&= \sum_i \big( (\omega_{\xi_2,\eta}\otimes\iota)\Delta(a_i)\xi_1 \big|
   \xi_i \big) \\
&= \sum_i \big( (1\otimes a_i)W(\xi_2\otimes\xi_1) \big|
   W(\eta\otimes\xi_i) \big).
\end{align*}
We will compute the image of $U^*$ by taking linear combinations of $\xi_1,
\xi_2$ and $\eta$.  In particular, as $W$ is unitary, we may replace
$W(\xi_2\otimes\xi_1)$ by $\xi_2\otimes\xi_1$ in the above expression.
It follows that the image of $U^*$ is the closed linear span of
vectors of the form
\[ \sum_i (\omega\otimes\iota)(W)^* a_i(\xi) \otimes f_i
\qquad (\omega\in L^1(\G), \xi\in L^2(\G)). \]
Now, $\{ (\omega\otimes\iota)(W)^* : \omega\in L^1(\G) \}$ is dense in
$C_0(\hat\G)$, and so the result follows, as
$V(\xi) = \sum_i a_i(\xi)\otimes f_i$.
As $C_0(\hat\G)$ contains a bounded approximate identity, clearly the image
of $U^*$ contains the image of $V$.  As $U^*U$ is the orthogonal projection
onto the image of $U^*$ (as $U^*$ is an isometry) it follows immediately that
$U^*U V = V$.
\end{proof}

\begin{proposition}\label{prop:three}
The corepresentation $U$ is unitary.
\end{proposition}
\begin{proof}
Suppose that $\sum_i \xi_i\otimes f_i\in L^2(\G)\otimes H$ is orthogonal to
the image of $U^*$.  Let $\xi\in L^2(\G)$ and $x'\in L^\infty(\G)'$, so
by Lemma~\ref{lem:image_ustar}, for any $\hat a\in C_0(\hat\G)$,
\[ 0 = \Big( (\hat a\otimes 1)Vx'\xi \Big| \sum_i \xi_i\otimes f_i \Big)
= \sum_i \big( \hat a a_i x'\xi \big| \xi_i \big)
= \sum_i \big( \hat a x' a_i \xi \big| \xi_i \big). \]
As $C_0(\hat\G)$ is $\sigma$-weakly dense in $L^\infty(\hat\G)$, it follows
that
\[ \sum_i \big( \hat x x' a_i \xi \big| \xi_i \big) = 0
\qquad (\xi\in L^2(\G), \hat x\in L^\infty(\hat\G), x'\in L^\infty(\G)')), \]
By Theorem~\ref{thm:gen_all_bh}, this implies that
\[ 0 = \sum_i \big( Ta_i\xi \big| \xi_i \big)
= \Big( (T\otimes 1)V\xi \Big| \sum_i \xi_i\otimes f_i \Big)
\qquad (\xi\in L^2(\G), T\in\mc B(L^2(\G))). \]
However, we know that $\{ (x\otimes 1)V(\xi) : x\in\mc B_0(L^2(\G)),
\xi\in L^2(\G) \}$ is linearly dense in $L^2(\G)\otimes H$.
Hence $\sum_i \xi_i\otimes f_i = 0$, and so $U^*$ has dense range, as required.
\end{proof}

\begin{remark}
We started with a completely positive left multiplier of $L^1(\hat\G)$;
let $L$ be the adjoint, a normal completely positive map on $L^\infty(\hat\G)$.
We immediately used Proposition~\ref{prop:one} to extend $L$ to a normal
completely positive map $\Phi$ on all of $\mc B(L^2(\G))$.  Remember
that the representation $\Phi(x) = V^*(x\otimes 1)V$ is unique (up to unitary
isomorphism), as this dilation is \emph{minimal}.  This is equivalent to the
non-degeneracy condition that $\{ (x\otimes 1)V\xi : x\in\mc B_0(L^2(\G)),
\xi\in L^2(\G) \}$ is linearly dense in $L^2(\G) \otimes H$.

As $\Phi$ extends $L$, we hence have a normal Stinespring representation of $L$,
as $L(\hat x) = V^*(\hat x\otimes 1)V$.  The previous results show that
$\{ (\hat x\otimes 1)V(\xi) : \hat x\in L^\infty(\hat\G), \xi\in L^2(\G) \}$
is also linearly dense in $L^2(\G) \otimes H$.  It follows that we also
have a \emph{minimal} Stinespring dilation of the original multiplier $L$.
\end{remark}

\subsection{Recovering the multiplier}

In the previous section, we showed how to construct a unitary corepresentation
of $\G$ from a completely positive left multiplier of $L^1(\hat\G)$.  We
now show how to recover the multiplier from the corepresentation.

\begin{proposition}\label{prop:recover}
There is $\alpha_0\in H$ such that
$U^*(\xi\otimes\alpha_0) = V(\xi) = \sum_i a_i(\xi)\otimes f_i$ for all
$\xi\in L^2(\G)$.
\end{proposition}
\begin{proof}
Let our left multiplier be represented by $a_0\in M(C_0(\G))$, so that
$(1\otimes a_0) \hat W = (\Phi\otimes\iota)(\hat W) =
\sum_i (a_i^*\otimes 1)\hat W(a_i\otimes 1)$.  Equivalently,
$\sum_i (1\otimes a_i^*)\Delta(a_i) = a_0\otimes 1$.
For $\omega\in L^1(\G)$ and $\xi,\eta\in L^2(\G)$, observe that
\begin{align*}
\Big( U^*\big(\xi\otimes\sum_i \ip{a_i}{\omega}f_i\big) \Big|
V(\eta) \big) &=
\sum_i \big( (\omega\otimes\iota)\Delta(a_i)\xi \big| a_i(\eta) \big) \\
&= \Big( (\omega\otimes\iota)\big( (1\otimes a_i^*)\Delta(a_i) \big) \xi
  \big| \eta \big)
= \ip{a_0}{\omega} (\xi|\eta).
\end{align*}
As this holds for all $\xi,\eta$, it follows that the map
$\sum_i \ip{a_i}{\omega}f_i \mapsto \ip{a_0}{\omega}$ is bounded,
and so the Riesz representation theorem for
Hilbert spaces provides $\alpha_0\in H$ such that
\[ \Big( \sum_i \ip{a_i}{\omega}f_i \Big| \alpha_0 \Big)
= \ip{a_0}{\omega} \qquad (\omega\in L^1(\G)). \]
By continuity,
\[ \big( U^*(\xi\otimes\alpha) \big| V(\eta) \big)
= (\xi\otimes\alpha|\eta\otimes\alpha_0) \qquad
(\xi,\eta\in L^2(\G), \alpha\in H), \]
that is, $UV(\eta) = \eta\otimes\alpha_0$ for all $\eta\in H$.
By Lemma~\ref{lem:image_ustar}, as $U^*UV=V$, it follows that
$V(\eta) = U^*UV(\eta) = U^*(\eta\otimes\alpha_0)$ as required.
\end{proof}

We now take slices of $U$ against this vector $\alpha_0$, and find that this
constructs our original multiplier, in the sense of
Proposition~\ref{prop:inv_corep_mult}.

\begin{theorem}\label{thm:mainthm}
Let $L_*$ be a completely positive left multiplier of $L^1(\hat\G)$.
There is a unitary corepresentation $U$ of $\G$ on $H$, such that $L=(L_*)^*$
is induced by $U$, using $\alpha_0\in H$.
\end{theorem}
\begin{proof}
Form $U$ as above and form $\alpha_0$ as in the previous proposition.
It is immediate that $a_i = (\iota\otimes\alpha_{\alpha_0,f_i})(U^*)$
for all $i$.  So the multiplier constructed by
Proposition~\ref{prop:inv_corep_mult} for $\alpha_0$ is, on $L^\infty(\hat\G)$,
the map $\hat x \mapsto a_i^* \hat x a_i$, which is just $L$, as required.
\end{proof}

\begin{theorem}\label{thm:mainthm2}
Let $\G$ be a locally compact quantum group.
There is an isometric, order preserving bijection between the completely
positive multipliers of $L^1(\hat\G)$ and $C_0^u(\hat\G)^*_+$.
\end{theorem}
\begin{proof}
Members of $C_0^u(\hat\G)^*_+$ induce completely positive left
multipliers of $L^1(\hat\G)$ in the sense discussed in
Section~\ref{sec:universal_link}.  Conversely, any completely positive
left multiplier comes from a unitary corepresentation, and this is
associated to a member of $C_0^u(\hat\G)^*_+$
by Proposition~\ref{prop:mults_from_uni_dual}.  That this procedure
gives a \emph{bijection} follows as $L^1(\hat\G)$ is an \emph{essential}
ideal in $C_0^u(\hat\G)^*$, see \cite[Proposition~8.3]{daws2}.

If $\mu\in C_0^u(\hat\G)^*_+$ is a state, then suppose that
 $C_0^u(\hat\G)\subseteq \mc B(H)$ is the universal representation, so
$\mu=\omega_{\alpha,\alpha}$ for some $\alpha\in H$.  Then $\mc W$ can be
identified with a member of $\mc B(L^2(\G)\otimes H)$, and
Proposition~\ref{prop:mults_from_uni_dual}
and Proposition~\ref{prop:inv_corep_mult} show that left multiplication by
$\mu$ induces the completely positive multiplier $L$, where in particular,
\[ L(1) = (\iota\otimes\omega_{\alpha,\alpha})(\mc W\mc W^*)
= 1 \ip{\mu}{1} = 1. \]
So $\|L\|=1$, and hence our bijection is an isometry.

Finally, if $\mu\leq\lambda$ in $C_0^u(\hat\G)^*_+$ then form the associated
completely positive multipliers $L_\mu$ and $L_\lambda$.  Let $L$ be the
multiplier formed from $\lambda-\mu$, so by uniqueness, $L=L_\lambda-L_\mu$.
As $L$ is completely positive, $L_\lambda \geq L_\mu$.  The converse is simply
a case of reversing the argument.  Thus our bijection is order preserving.
\end{proof}

We remark that it is completely obvious from these results that for any
completely \emph{positive} left multiplier $L$, there is a completely
positive right multiplier $L'$ such that $(L,L')$ forms a double multiplier
(simply let $L'$ be induced by right multiplication by the element of
$C_0^u(\hat\G)^*$ associated to $L$).  It seems to be unknown if a similar
result holds for completely bounded multipliers.

\section{Representing elements for completely bounded multipliers}

Notice that while Proposition~\ref{prop:mults_rep} shows that all completely
bounded multipliers are ``represented'', we didn't use this fact until 
Proposition~\ref{prop:recover}.  Here we show how to use
the representing element more directly.

Recall, from \cite{kus} for example, that $L^1(\G)$ contains a dense
$*$-subalgebra $L^1_\sharp(\G)$; we define $\omega\in L^1_\sharp(\G)$ if
and only if there is $\tau\in L^1(\G)$ with $\ip{x}{\tau} =
\overline{ \ip{S(x)^*}{\omega} }$ for all $x\in D(S)$, and in this case,
denote $\omega^\sharp = \tau$.  We note that the elementary
properties of $L^1_\sharp(\G)$ can be developed \emph{mutatis mutandis} for
$\G$ coming from manageable multiplicative unitaries.

Recall that the scaling group $(\tau_t)$ is implemented as $\tau_t(x) =
P^{it} x P^{-it}$, where $P$ is a certain positive injective operator.
As $R$ and $\tau_t$ commute for all $t$, and $S=R\tau_{-i/2}$, it follows
that $R$ leaves $D(S)$ invariant, and $RS=SR$.  It is then easy to see
that $R_*$ leaves $L^1_\sharp(\G)$ invariant, and $R_*(\omega^\sharp)
= R_*(\omega)^\sharp$ for $\omega\in L^1_\sharp(\G)$.
Given $\beta\in D(P^{-1/2})$ and $\xi\in D(P^{1/2})$, we have that for
$x\in D(S) = D(\tau_{-i/2})$,
\[ \ip{x}{\omega_{P^{-1/2} \beta, P^{1/2} \xi}}
= \big( P^{1/2} x P^{-1/2} \beta \big| \xi \big)
= \ip{\tau_{-i/2}(x)}{\omega_{\beta,\xi}}
= \ip{S(R(x))}{\omega_{\xi,\beta}^*}
= \ip{x}{R_*(\omega_{\xi,\beta}^\sharp)}, \]
and so $\omega_{\xi,\beta}\in L^1_\sharp(\G)$ with
$\omega^\sharp_{\xi,\beta} = R_*(\omega_{P^{-1/2} \beta, P^{1/2} \xi})$.

\begin{proposition}\label{prop:rep_gives_mult}
Let $L$ be a completely bounded left multiplier of $L^1(\hat\G)$, represented
by $a_0\in M(C_0(\G))$.  For $\xi,\eta\in D(P^{1/2})$ and
$\alpha,\beta\in D(P^{-1/2})$, we have that
\[ \big( \Phi(\theta_{\xi,\eta})\alpha \big| \beta \big)
= \ip{\Delta(a_0)}{\omega_{\alpha,\eta} \otimes
\omega_{\xi,\beta}^{\sharp*}}. \]
\end{proposition}
\begin{proof}
Let $\xi_0\in L^2(\G)$ be a unit vector, let $(e_i)$ be an orthonormal basis
for $L^2(\G)$, and let $W(\alpha\otimes\xi_0) = \sum_i \alpha_i \otimes e_i$
and $W(\beta\otimes\xi_0) = \sum_i \beta'_i\otimes e_i$.
For $\epsilon>0$, we can find a family $(\beta_i)$ in $D(P^{-1/2})$ with
\[ \Big\| W(\beta\otimes\xi_0) - \sum_i \beta_i\otimes e_i \Big\|<\epsilon. \]
Using Proposition~\ref{prop:one}, and that $\hat W = \Sigma W^* \Sigma$,
we see that
\begin{align*} \big( \Phi(\theta_{\xi,\eta})\alpha \big| \beta \big)
&= \big( (\iota\otimes L)(W(\theta_{\xi,\eta}\otimes 1)W^*) W(\alpha\otimes\xi_0)
   \big| W(\beta\otimes\xi_0) \big) \\
&= \sum_{i,j} \big( (\omega_{\alpha_i,\beta_j'}\otimes\iota)
   (\iota\otimes L)(W(\theta_{\xi,\eta}\otimes 1)W^*) e_i \big| e_j \big) \\
&= \sum_{i,j} \big( L( (\omega_{\xi,\beta_j'}\otimes\iota)(W)
   (\omega_{\alpha_i,\eta}\otimes\iota)(W^*) ) e_i \big| e_j \big).
\end{align*}
A similar calculation establishes that if
\[ x = \sum_{i,j} \big( L( (\omega_{\xi,\beta_j}\otimes\iota)(W)
(\omega_{\alpha_i,\eta}\otimes\iota)(W^*) ) e_i \big| e_j \big), \]
then
\[ \big| \big( \Phi(\theta_{\xi,\eta})\alpha \big| \beta \big) - x \big|
< \epsilon \|L\|_{cb} \|\alpha\| \|\xi\| \|\eta\|. \]
That is, we may replace $(\beta_j')$ by $(\beta_j)$, at the cost of
a small error term.

As $(\omega\otimes\iota)(W)^* = (\omega^\sharp\otimes\iota)(W)$ for
$\omega\in L^1_\sharp(\G)$, we see that $(\omega_{\xi,\beta_j}\otimes\iota)
(W) = (\omega_{\xi,\beta_j}^\sharp\otimes\iota)(W)^*
= (\omega_{\xi,\beta_j}^{\sharp *}\otimes\iota)(W^*)$.  This makes sense,
as $\beta_j\in D(P^{-1/2})$ and $\xi\in D(P^{1/2})$.  Thus
\[ (\omega_{\xi,\beta_j}\otimes\iota)(W)
   (\omega_{\alpha_i,\eta}\otimes\iota)(W^*)
= (\omega_{\xi,\beta_j}^{\sharp *}\otimes\iota)(W^*)
(\omega_{\alpha_i,\eta}\otimes\iota)(W^*)
= (\omega_{\alpha_i,\eta}\omega_{\xi,\beta_j}^{\sharp *}\otimes\iota)(W^*). \]
Recall that $(\iota\otimes L)(W^*) = (a_0\otimes 1)W^*$, and that
$(\Delta\otimes\iota)(W^*) = W^*_{23} W^*_{13}$, and so
\begin{align*} x
&= \sum_{i,j} \big( (\omega_{\alpha_i,\eta}
   \omega_{\xi,\beta_j}^{\sharp *}\otimes\iota)
   ( (\iota\otimes L)(W^*) ) e_i \big| e_j \big) \\
&= \sum_{i,j} \ip{(\Delta(a_0)\otimes 1)W^*_{23} W^*_{13}}{\omega_{\alpha_i,\eta}
\otimes\omega_{\xi,\beta_j}^{\sharp *}\otimes \omega_{e_i,e_j}} \\
&= \sum_j \ip{(\Delta(a_0)\otimes 1)W^*_{23}}{\omega_{\alpha,\eta}
\otimes \omega_{\xi,\beta_j}^{\sharp *} \otimes \omega_{\xi_0,e_j}}.
\end{align*}

Let $a\in D(S)^*$, so that
\begin{align*}
\sum_j & \ip{(a\otimes 1)W^*}{\omega_{\xi,\beta_j}^{\sharp *}
   \otimes \omega_{\xi_0,e_j}}
= \sum_j \ip{a S((\iota\otimes\omega_{\xi_0,e_j})(W))}{
   \omega_{\xi,\beta_j}^{\sharp *}} \\
&= \sum_j \overline{ \ip{S((\iota\otimes\omega_{\xi_0,e_j})(W))^* a^*}{
   \omega_{\xi,\beta_j}^\sharp} }
= \sum_j \ip{(\iota\otimes\omega_{\xi_0,e_j})(W) S(a^*)^*}{
   \omega_{\xi,\beta_j}} \\
&= \sum_j \ip{W (S(a^*)^*\otimes 1)}{\omega_{\xi,\beta_j}
   \otimes \omega_{\xi_0,e_j}}
= \Big( W(S(a^*)^*\otimes 1) (\xi\otimes\xi_0) \big|
   \sum_j \beta_j\otimes e_j \Big). \end{align*}
By comparison,
\begin{align*} \Big( W(S(a^*)^*\otimes 1) (\xi\otimes\xi_0) \big|
   \sum_j \beta'_j\otimes e_j \Big)
&= \big( (S(a^*)^*\otimes 1)(\xi\otimes\xi_0)\big|
   W^*W(\beta\otimes\xi_0) \big) \\
&= \ip{S(a^*)^*}{\omega_{\xi,\beta}}
= \overline{ \ip{a^*}{\omega_{\xi,\beta}^\sharp} }
= \ip{a}{\omega_{\xi,\beta}^{\sharp*}}.
\end{align*}
If it so happens that $a = (\omega_{\alpha,\eta}\otimes\iota)\Delta(a_0)$
is in $D(S)^*$, then we have
\[ \big| x - \ip{a}{\omega_{\xi,\beta}^{\sharp*}} \big|
\leq \epsilon \|\xi\| \|S(a^*)\|. \]
However, observe that for this choice of $a$,
\[ \ip{a}{\omega_{\xi,\beta}^{\sharp*}}
= \ip{\Delta(a_0)}{\omega_{\alpha,\eta} \otimes
\omega_{\xi,\beta}^{\sharp*}}, \]
and so as $\epsilon>0$, this gives the required result.

So it remains to show that $a=(\omega_{\alpha,\eta}\otimes\iota)\Delta(a_0)
\in D(S)^*$.  By \cite[Theorem~5.9]{daws1}, we know that $a_0\in D(S)^*$,
and by hypothesis, $\omega_{\eta,\alpha}\in L^1_\sharp(\G)$.  Thus, for
$\omega\in L^1_\sharp(\G)$,
\[ \ip{a^*}{\omega^\sharp}
= \ip{(\omega_{\eta,\alpha}\otimes\iota)\Delta(a_0^*)}{\omega^\sharp}
= \ip{a_0^*}{\omega_{\eta,\alpha} \omega^\sharp}
= \overline{ \ip{S(a_0^*)^*}{\omega \omega_{\eta,\alpha}^\sharp} }
= \overline{ \ip{(\iota\otimes\omega_{\eta,\alpha}^\sharp)\Delta(S(a_0^*)^*)}
{\omega} }. \]
This is enough to show that $a^*\in D(S)$ with $S(a^*)
= (\iota\otimes\omega_{\eta,\alpha}^\sharp)\Delta(S(a_0^*)^*)$, given that
$S$ is a $\sigma$-weakly closed operator; for details see for
example \cite[Appendix~A]{bds}.
\end{proof}

\subsection{Weak$^*$-continuity of the Junge, Neufang, Ruan representation}

As explained in Section~\ref{sec:mult} above, \cite{jnr} shows that for
a locally compact quantum group $\G$, there is a bijection between the
completely bounded left multipliers of $L^1(\hat\G)$, say
$M_{cb}^l(L^1(\hat\G))$, and
$\mc{CB}^{\sigma, L^\infty(\hat\G)}_{L^\infty(\G)'}(\mc B(L^2(\G)))$.
In \cite{hnr2}, it is shown that this map is weak$^*$-weak$^*$ continuous,
at least when $\hat\G$ has the \emph{left co-AP} property, see
\cite[Corollary~4.10]{hnr2} (and \cite[Theorem~4.7]{hnr2} for the version
for right multipliers).  In this final section of the paper, we apply the
result of the previous section to show that
this weak$^*$-continuity result is true for all $\G$.

Firstly, we recall from \cite{hnr2} the proof that $M_{cb}^l(L^1(\hat\G))$
is a dual space.  Proposition~\ref{prop:mults_rep} shows that we have a
map $\Lambda:M_{cb}^l(L^1(\hat\G)) \rightarrow L^\infty(\G)$ (actually,
this maps into $M(C_0(\G))$, but this is unimportant here) which satisfies
\[ (L\otimes\iota)(W)W^* = 1\otimes\Lambda(L_*), \qquad
\Lambda(L_*) \hat\lambda(\hat\omega) =
\hat\lambda\big( L_*(\hat\omega) \big)
\qquad (\hat\omega \in L^1(\hat\G)). \]
It follows that $\Lambda$ is a contractive algebra homomorphism.
Then \cite[Proposition~3.4]{hnr2} shows that if we denote by $X$ the image
of $\Lambda$, equipped with the norm coming from $M_{cb}^l(L^1(\hat\G))$,
then the closed unit ball of $X$ is weak$^*$-closed in $L^\infty(\G)$.  Indeed,
giving $M_{cb}^l(L^1(\hat\G))$ its canonical operator space structure, the
closed unit ball of $M_n(X)$ is weak$^*$-closed in $M_n(L^\infty(\G))$.
Using this, \cite[Theorem~3.5]{hnr2} shows that if we let
$Q^l_{cb}(L^1(\hat\G))$ be the closure in $M_{cb}^l(L^1(\hat\G))^*$
of the image of $L^1(\G)$ under the adjoint of $\Lambda$, then 
$Q^l_{cb}(L^1(\hat\G))^*$ is completely isometrically isomorphic to
$M_{cb}^l(L^1(\hat\G))$.  Thus we get a weak$^*$-topology on
$M_{cb}^l(L^1(\hat\G))$.

In \cite[Section~8]{daws2} we independently gave an analogous
construction of a weak$^*$-topology on the space of double multipliers.
In fact, the first part of the proof of \cite[Proposition~8.11]{daws2}
already works for merely left multipliers, and then one can apply the
abstract result which is \cite[Proposition~8.12]{daws2} to construct
$Q^l_{cb}(L^1(\hat\G))$.  In \cite{daws2} we found a very ``Banach
algebraic'' way to construct preduals for double multiplier algebras
(see \cite[Theorem~7.7]{daws2} for example), but it seems that at several
crucial points, it really is necessary to work with double multipliers.
It would be interesting to know how to adapt these ideas to one-sided
multipliers.

For us, the important point is that if $(L_\alpha)$ is a bounded net
in $M_{cb}^l(L^1(\hat\G))$, then $(L_\alpha)$ is weak$^*$-null with
respect to $Q^l_{cb}(L^1(\hat\G))$ if and only if $(\Lambda(L_\alpha))$
is weak$^*$-null in $L^\infty(\G)$.

We next consider the space
$\mc{CB}^{\sigma, L^\infty(\hat\G)}_{L^\infty(\G)'}(\mc B(L^2(\G)))$.
Firstly, we consider the larger space
$\mc{CB}^{\sigma}(\mc B(L^2(\G)))$ which can be identified with
$\mc{CB}(\mc B_0(L^2(\G)), \mc B(L^2(\G)))$.  This in turn is the dual
space of $\mc B_0(L^2(\G)) \proten \mc B(L^2(\G))_*$, the operator space
projective tensor product of the compact operators $\mc B_0(L^2(\G))$ with the
trace-class operators $\mc B(L^2(\G))_*$.  By restriction, we have a
weak$^*$-topology on 
$\mc{CB}^{\sigma, L^\infty(\hat\G)}_{L^\infty(\G)'}(\mc B(L^2(\G)))$.
Again, for us the important point is that a bounded net $(\Phi_\alpha)$
in $\mc{CB}^{\sigma, L^\infty(\hat\G)}_{L^\infty(\G)'}(\mc B(L^2(\G)))$
is weak$^*$-null if and only if $(\Phi_\alpha(\theta))$ is a weak$^*$-null
net in $\mc B(L^2(\G))$, for each $\theta\in\mc B_0(L^2(\G))$.
All this is explained in \cite[Section~4]{hnr2} and the references therein.

The following improves \cite[Theorem~4.7]{hnr2} (which is stated for
right multipliers) in that we need make no approximation property type
assumptions.

\begin{theorem}\label{thm:weakstar}
For any $\G$, the map $M_{cb}^l(L^1(\hat\G)) \rightarrow
\mc{CB}^{\sigma, L^\infty(\hat\G)}_{L^\infty(\G)'}(\mc B(L^2(\G)))$ is
weak$^*$-weak$^*$-continuous.  If $\G$ is a locally compact quantum
group, this correspondence is a weak$^*$-weak$^*$-continuous homeomorphism.
\end{theorem}
\begin{proof}
Denote by $\phi$ the map $M_{cb}^l(L^1(\hat\G)) \rightarrow
\mc{CB}^{\sigma, L^\infty(\hat\G)}_{L^\infty(\G)'}(\mc B(L^2(\G)))$.
To show that $\phi$ is weak$^*$-continuous,
it suffices to show that if $(L_i)$ is a bounded, weak$^*$-null
net in $M_{cb}^l(L^1(\hat\G))$, then the corresponding bounded net,
say $(\Phi_i)$, in $\mc{CB}^\sigma(\mc B(L^2(\G)))$ is weak$^*$-null.
When $\G$ is a locally compact quantum group, we know from \cite{jnr} that
$\phi$ is a completely isometric isomorphism, and then if $\phi$ is
weak$^*$-continuous, it is automatically a weak$^*$-weak$^*$-continuous 
homeomorphism.  This is perhaps not well-known (in the operator space
setting) but see \cite[Lemma~10.1]{daws2} for example.

We fix a bounded weak$^*$-null net $(L_i)$ of left multipliers, with
corresponding net $(\Phi_i)$.  For each $i$ let $L_i$
be represented by $a_i\in L^\infty(\G)$.  That $(L_i)$ is weak$^*$-null
means that $(a_i)$ is weak$^*$-null in $L^\infty(\G)$.  As explained above,
as $\Lambda$ is a contraction, $(a_i)$ is also a bounded net.
By Proposition~\ref{prop:rep_gives_mult}, we have that
\[ (\Phi_i(\theta_{\xi,\eta})\alpha|\beta)
= \ip{a_i}{\omega_{\alpha,\eta} \omega_{\xi,\beta}^{\sharp *}}
\qquad (\xi,\eta\in D(P^{1/2}), \alpha,\beta\in D(P^{-1/2})). \]
As $D(P^{1/2})$ and $D(P^{-1/2})$ are dense in $L^2(\G)$,
we immediately see that
\[ \lim_i \ip{\Phi_i(\theta)}{\omega} = 0 \]
for a dense collection of $\theta\in\mc B_0(L^2(\G))$ and
$\omega\in\mc B(L^2(\G))_*$.  As $(\Phi_i)$ is a bounded net, this is
enough to show that $(\Phi_i)$ is weak$^*$-null, as required.
\end{proof}

As remarked on in the proof of \cite[Theorem~4.7]{hnr2}, we can equivalently
state this result in terms of the Haagerup tensor product
(see \cite[Chapter~9]{er}).  We have a completely isometric isomorphism
\[ \mc B_0(L^2(\G)) \proten \mc B(L^2(\G))_* \rightarrow
\mc B(L^2(\G))_* \overset{h}{\otimes} \mc B(L^2(\G))_*; \qquad
\theta_{\xi,\eta} \otimes \omega_{\alpha,\beta} \mapsto
\omega_{\xi,\beta} \otimes \omega_{\alpha,\eta}. \]
See also the discussion after \cite[Remark~4.6]{hnr2}.
The adjoint gives a normal completely isometric isomorphism
\[ \mc B(L^2(\G)) \overset{eh}{\otimes} \mc B(L^2(\G))
\rightarrow \mc{CB}^\sigma(\mc B(L^2(\G))); \quad
x\otimes y \mapsto T_{x,y}. \]
Here we use the extended (or weak$^*$) Haagerup tensor product,
see \cite{bsm, er2}, and $T_{x,y}$ is the operator $z\mapsto xzy$.
This isomorphism restricts to an
isomorphism between $L^\infty(\G) \overset{eh}{\otimes} L^\infty(\G)$
and $\mc{CB}^\sigma_{L^\infty(\G)'}(\mc B(L^2(\G)))$, and the predual
of $L^\infty(\G) \overset{eh}{\otimes} L^\infty(\G)$ is
$L^1(\G) \overset{h}{\otimes} L^1(\G)$.

We can hence restate Theorem~\ref{thm:weakstar} as saying that there is
a completely bounded map
\[ \phi_*: L^1(\G) \overset{h}{\otimes} L^1(\G) \rightarrow
Q^l_{cb}(L^1(\hat\G)), \]
the adjoint of which is our map $M_{cb}^l(L^1(\hat\G)) \rightarrow
\mc{CB}^\sigma_{L^\infty(\G)'}(\mc B(L^2(\G)))$.

\small
\noindent
Matthew Daws\\
School of Mathematics,\\
University of Leeds,\\
LEEDS LS2 9JT\\
United Kingdom\\
Email: \texttt{matt.daws@cantab.net}

\end{document}